\documentclass[11pt]{amsart}
\usepackage[T1]{fontenc}
\usepackage{lmodern,microtype,amssymb,mathtools,mathrsfs}
\usepackage{enumitem,needspace}
\usepackage[dvipsnames]{xcolor}
\usepackage[colorlinks=true,linkcolor=MidnightBlue,citecolor=MidnightBlue,urlcolor=MidnightBlue]{hyperref}
\hypersetup{pdftitle={Holomorphic integer-graded vertex operator superalgebras are purely even},pdfauthor={Chongying Dong, Kaiwen Sun and Haowu Wang}}
\setlist{itemsep=3pt,topsep=6pt}
\newtheorem{theorem}{Theorem}[section]
\newtheorem{proposition}[theorem]{Proposition}
\newtheorem{lemma}[theorem]{Lemma}
\newtheorem{corollary}[theorem]{Corollary}
\theoremstyle{definition}
\newtheorem{definition}[theorem]{Definition}
\theoremstyle{remark}
\newtheorem{remark}[theorem]{Remark}

\newcommand{\C}{\mathbb C}
\newcommand{\Z}{\mathbb Z}
\newcommand{\Q}{\mathbb Q}
\newcommand{\N}{\mathbb Z_{\geq0}}
\newcommand{\one}{\mathbf1}
\newcommand{\Cat}{\mathcal C}
\newcommand{\Ind}{\mathcal I}
\newcommand{\Hom}{\operatorname{Hom}}
\newcommand{\End}{\operatorname{End}}
\newcommand{\Rep}{\operatorname{Rep}}

\newcommand{\Tr}{\operatorname{Tr}}
\newcommand{\Span}{\operatorname{span}}
\newcommand{\wt}{\operatorname{wt}}
\newcommand{\id}{\operatorname{id}}
\newcommand{\qdim}{\operatorname{qdim}}
\newcommand{\M}{\mathsf M}
\newcommand{\sVec}{\mathrm{sVec}}
\numberwithin{equation}{section}
\title[Holomorphic integer-graded vertex superalgebras]{Holomorphic integer-graded vertex operator superalgebras are purely even}

\author[C. Dong]{Chongying Dong}
\address{Department of Mathematics, University of California,
Santa Cruz, CA 95064, USA}
\email{dong@ucsc.edu}

\author[K. Sun]{Kaiwen Sun}
\address{School of Mathematical Sciences, University of Science and
Technology of China, Hefei 230026, Anhui, China}
\email{kwsun@ustc.edu.cn}

\author[H. Wang]{Haowu Wang}
\address{School of Mathematics and Statistics, Wuhan University,
Wuhan 430072, Hubei, China}
\email{haowu.wangmath@whu.edu.cn}

\date{}
\begin{document}
\begin{abstract}
We prove that a simple, self-contragredient, rational, $C_2$-cofinite, integer-graded holomorphic vertex operator superalgebra of CFT-type has zero odd part. We first show that a hypothetical nonzero odd part would be a simple current of categorical dimension $-1$ over a strongly rational pointed VOA with four simple modules. We then prove that every nonvacuum simple current of a strongly rational VOA has strictly positive lowest conformal weight. This leads to a contradiction with Dong--Jiao--Xu's theorem on the positivity of categorical dimension. 
\end{abstract}

\maketitle
%\tableofcontents

\section{Introduction}
\label{sec:intro}
\renewcommand{\thetheorem}{\arabic{theorem}}

Holomorphic vertex operator algebras are distinguished for that their only irreducible module is themselves. Their representation theory is particularly simple, while their internal structure can be remarkably rich. The moonshine module constructed by Frenkel, Lepowsky and Meurman \cite{FLM} is the most famous example. At central charge $24$, Schellekens \cite{Schellekens} obtained a list of $71$ possible weight-one Lie algebras for holomorphic VOAs.
The construction and uniqueness of the holomorphic VOAs with nonzero weight-one space have since been established, see
\cite{EMS,BLS,MSClassification} and the references therein. The uniqueness of moonshine module remains a major open question. 

For vertex operator superalgebras, one must distinguish two conformal grading conventions. In the usual half-integer grading, even vectors have integer weights and odd vectors have weights in $\frac12+\Z$. Supersymmtric minimal models and free fermion are familiar examples. This grading corresponds to the chiral algebras for what physicists usually call two-dimensional fermionic CFTs. Creutzig, Duncan and Riedler \cite{CDR} studied self-dual vertex superalgebras of small central charge, including their classification at central charge $12$.
H\"ohn and M\"oller \cite{HM} classified the self-dual superalgebras of central charge at most $24$ under regularity hypotheses including positive lowest weights for nonvacuum untwisted modules and all parity-twisted modules. 
At central charge $24$ their list contains $969$ superalgebras with nonzero odd part, subject to uniqueness of the shorter moonshine module. 

In an integer-graded vertex operator superalgebra, parity and conformal weight are independent: an odd vector may have integer weight. Affine vertex superalgebras provide many examples, where the weight-one space is now a Lie superalgebra rather than necessarily a Lie algebra. 
Van Ekeren and Rodr\'iguez Morales \cite{EM} investigated holomorphic vertex operator superalgebras with integer grading. They proved that those of central charge $8$ or $16$ are purely even. More precisely, the possibilities are the lattice VOA
$V_{E_8}$ at central charge $8$, and $V_{E_8}\otimes V_{E_8}$ or
$V_{D_{16}^{+}}$ at central charge $16$. 

The analysis in \cite{EM} at central charge $24$ begins with the following superalgebra analogue of Schellekens' first trace identity. Let $V$ be a holomorphic integer-graded VOSA of CFT-type satisfying the finiteness and duality assumptions of \cite{EM}. Van Ekeren and Rodr\'iguez Morales show that in the case of $\operatorname{sdim}V_1\neq 24$, $V_1$ is zero
or a direct sum of simple Lie superalgebras, each either purely even or isomorphic to $\mathfrak{osp}(1|2n)$. For the nonzero simple summands, their positive integral levels $k_i$ and dual Coxeter numbers $h_i^\vee$ satisfy
\begin{equation}\label{eq:level-identity}
 \frac{h_i^\vee}{k_i}
   =\frac{\operatorname{sdim}V_1-24}{24}.
\end{equation}
Their original enumeration reports $1332$ candidate semisimple weight-one Lie superalgebras satisfying these restrictions.\footnote{The count of candidates with nonzero odd part is corrected to $1106$ in \cite{SWW}.}
This is a list of candidates, not a list of existing holomorphic vertex superalgebras. Here the speciality of $\mathfrak{osp}(1|2n)$ is due to the fact that simple affine VOSA associated with a basic simple Lie superalgebra with nonzero odd part is rational and $C_2$-cofinite if and only if the underlying algebra is Lie superalgebra of type $\mathfrak{osp}(1|2n)$ and the level is a positive integer \cite{GK,AL,CLosp}.

The existence of a holomorphic integer-graded VOSA with nonzero odd part was left open in \cite{EM}. In this paper, we prove that no such VOSA exists under the standard simplicity, rationality, cofiniteness and duality assumptions.
Our argument actually applies to every central charge. 
Throughout this paper, rationality means complete reducibility of admissible untwisted modules. Holomorphicity means that the adjoint module is the only irreducible admissible untwisted module up to parity reversal.
Self-contragredience means the existence of a nondegenerate even
invariant bilinear form. The definitions will be recalled in
Section~\ref{sec:prelim}.

\begin{theorem}\label{thm:general}
Every simple, self-contragredient, rational, $C_2$-cofinite, holomorphic integer-graded vertex operator superalgebra of CFT-type is purely even.
\end{theorem}

\begin{corollary}\label{cor:24}
A simple, self-contragredient, rational, $C_2$-cofinite, holomorphic
integer-graded VOSA of CFT-type and central charge $24$ has zero odd part.
\end{corollary}

The Corollary~\ref{cor:24} provides a VOA counterpart to a number theory theorem on automorphic forms \cite[Theorem~11.5]{SWW}. 
The holomorphic VOAs with central charge 24 have been known to closely relate to certain distinguished automorphic forms, called automorphic products or Borcherds products \cite{Borcherds1995,Borcherds1998}. 
It was proved in \cite{DSW23,SWW23} that the holomorphic VOAs with central charge 24 have one-to-one correspondence to the anti-symmetric holomorphic automorphic products of singular weight on $2U\oplus L$ at one-dimensional cusp, where $U$ is a hyperbolic plane and $L$ is an even positive-definite lattice of positive rank. However, it was an open question whether super generalizations can be allowed on each side. Recently, Sun, Wang and Williams proved in \cite{SWW} that an anti-symmetric holomorphic automorphic product of singular weight on
$2U\oplus L$ has no odd real roots in its associated
Borcherds--Kac--Moody superalgebra. This means all such anti-symmetric automorphic products are exactly those already given in \cite{DSW23,SWW23}.
The anti-symmetric automorphic product satisfies the same relation
\eqref{eq:level-identity}, see \cite[Theorem~1.6]{SWW}. In \cite{SWW}, more than one thousand candidates with odd part were ruled out by heavy calculations and a large number of case-by-case studies. 
Our argument in this paper establishes the vertex-algebraic nonexistence result directly, without assuming an automorphic-product construction. 

A second main result, used in the proof and of independent interest, concerns the conformal weights of simple currents. Recall that an ordinary VOA is \emph{strongly rational} if it is simple, self-contragredient, rational, $C_2$-cofinite and of CFT-type. A \emph{simple current} is an irreducible module invertible under fusion.

\begin{theorem}\label{thm:positive-currents}
Let $A$ be a strongly rational VOA. Every simple current $J\not\cong A$
has strictly positive lowest conformal weight.
\end{theorem}

The conclusion concerns the actual lowest weight $h_J$, not its class
modulo $\Z$. It follows from an inequality for four-point functions.
For every irreducible $A$-module $M$, we prove
\begin{equation}\label{eq:intro-weight-inequality}
 2h_J+2h_M-h_{J\boxtimes M}-h_{J^{-1}\boxtimes M}\in\N.
\end{equation}
Summing this inequality along a finite orbit of fusion by $J$ gives
$h_J\geq0$. If $h_J=0$, equality holds throughout every orbit. Balancing identity of modular tensor category then makes $J$ in M\"uger center, and nondegeneracy of $\Rep(A)$ implies
$J\cong A$. This proves strict positivity for a nonvacuum simple current.

We explain how these results fit together. Suppose that $V$ satisfies
the hypotheses of Theorem~\ref{thm:general}, and write
\begin{equation}\label{eq:UF}
 U=V_{\bar0},\qquad F=V_{\bar1},\qquad V=U\oplus F.
\end{equation}
Assuming $F\neq0$, we first establish the following structural fact.

\begin{proposition}\label{prop:descent-intro}
The even subalgebra $U$ is strongly rational. It has exactly four
irreducible modules, all simple currents, with fusion group
$\Z_2\times\Z_2$.
\end{proposition}

The proof begins with the $C_2$-cofiniteness of $U$, using the parity version of Miyamoto's cyclic-orbifold argument \cite{MiyCyclic}. In the resulting braided category of generalized $U$-modules, the odd summand $F$ is invertible and has self-braiding $-1$. Monodromy with $F$ splits the category into two components. Induction to $V$ and holomorphicity show that the neutral component is semisimple, with simple objects $U$ and $F$.
Contragredient duality then gives an invertible object in the other
component and proves its semisimplicity. Thus rationality of $U$ is a consequence of the hypotheses on $V$.

The integer grading gives twist $\theta_F=1$, while the odd-odd supercommutativity gives self-braiding $c_{F,F}=-1$. The spin--statistics identity \cite[Corollary~2.8]{CKL} then implies categorical dimension $d_F=-1$. On the other hand, all nonvacuum simple $U$-modules have positive lowest conformal weight by Theorem~\ref{thm:positive-currents}.
The quantum-dimension formula of Dong--Jiao--Xu \cite{DJX}, together
with Huang's modularity theorem \cite{Huang}, identifies $d_F$ with a
positive character-ratio limit. This is the required contradiction.

\begin{remark}
The holomorphicity assumption is essential. For example, the
integer-graded affine VOSA $L_1(\mathfrak{osp}(1|2n))$ has nonzero odd
part, is strongly rational, and has $n+1$ irreducible untwisted modules
up to parity reversal, see \cite{CLosp,CGLosp}. Thus integer grading alone
does not obstruct rational vertex superalgebras.
\end{remark}

The paper is organized as follows. 
Section~\ref{sec:prelim} recalls the required definitions and tensor
category results. Sections~\ref{sec:c2} and \ref{sec:descent} establish
Proposition~\ref{prop:descent-intro}. We prove
Theorem~\ref{thm:positive-currents} in Section~\ref{sec:positivity},
and conclude the proof of Theorem~\ref{thm:general} in
Section~\ref{sec:characters}. As byproducts of Theorem~\ref{thm:positive-currents}, we prove a sharp lower bound for the conformal weight of nonvacuum simple current in Corollary \ref{cor:bound}, and prove M\"oller's positivity conjecture \cite[Conjecture~1.1]{Moeller} in Corollary \ref{cor:moeller}.

\renewcommand{\thetheorem}{\thesection.\arabic{theorem}}
\section{Vertex operator superalgebras and their modules}
\label{sec:prelim}

We follow the integer-grading conventions of
\cite[Section~4]{EM}. In particular, conformal weight does not determine
parity.

\subsection{Vertex operator superalgebras}

A vector superspace is a vector space $V=V_{\bar0}\oplus V_{\bar1}$.
For a homogeneous vector $a$, its parity is denoted by $|a|\in\{0,1\}$.
A vertex superalgebra consists of such a space, an even vacuum vector
$\one$, an even translation operator $D$, and a parity-preserving
state-field correspondence
$$
 Y(a,z)=\sum_{m\in\Z}a_mz^{-m-1},
 \qquad Y(a,z)b\in V((z)).
$$
The vacuum and translation axioms are
\begin{align*}
 Y(\one,z)&=\id_V,& Y(a,z)\one&\in a+zV[[z]],\\
 D\one&=0,& [D,Y(a,z)]&=\partial_zY(a,z).
\end{align*}
For homogeneous $a,b$, the super Jacobi identity is
\begin{align}
 &z_0^{-1}\delta\left(\frac{z_1-z_2}{z_0}\right)
       Y(a,z_1)Y(b,z_2)\notag\\
 &\quad-(-1)^{|a||b|}
 z_0^{-1}\delta\left(\frac{z_2-z_1}{-z_0}\right)
       Y(b,z_2)Y(a,z_1)\notag\\
 &\hspace{2cm}
 =z_2^{-1}\delta\left(\frac{z_1-z_0}{z_2}\right)
       Y(Y(a,z_0)b,z_2),\label{eq:super-jacobi}
\end{align}
where $\delta(z)=\sum_{m\in\Z}z^m$ and binomial expressions have their usual formal expansions.

A conformal vector is an even vector $\omega$ such that
$$
 Y(\omega,z)=\sum_{n\in\Z}L(n)z^{-n-2},\qquad L(-1)=D,
$$
and
$$
 [L(m),L(n)]
 =(m-n)L(m+n)+\frac{m^3-m}{12}\delta_{m+n,0}\,c\,\id_V
$$
for a central charge $c$. An \emph{integer-graded vertex
operator superalgebra} is a conformal vertex superalgebra for which
$L(0)$ is semisimple and
$$
 V=\bigoplus_{n\in\Z}V_n,\qquad V_n=\{v:L(0)v=nv\}.
$$
We require $\dim V_n<\infty$ and $V_n=0$ for $n\ll0$.
Each $V_n$ has its own parity decomposition. An ordinary VOA is the
special case $V_{\bar1}=0$. We say that $V$ is of \emph{CFT-type} if $V_n=0$ for $n<0$ and $V_0=\C\one$. Ideals are understood to be
parity-graded, and simplicity means the absence of proper nonzero ideals.

For vectors homogeneous in conformal weight and parity, the weight
formula and skew symmetry are
\begin{equation}\label{eq:weight-mode}
 \wt(a_mb)=\wt(a)+\wt(b)-m-1,
\end{equation}
whenever $a_mb\neq0$, and
\begin{equation}\label{eq:skew}
 Y(a,z)b=(-1)^{|a||b|}e^{zD}Y(b,-z)a.
\end{equation}
The corresponding supercommutator formula is
$$
 [a_m,b_n]_{\mathrm{super}}
 =\sum_{j\geq0}\binom mj(a_jb)_{m+n-j}.
$$
For a VOSA of CFT-type, $V_1$ is a finite-dimensional Lie superalgebra
with bracket $[a,b]=a_0b$. If $T$ is an even endomorphism
of a finite-dimensional superspace $W$, its supertrace is
$\operatorname{str}_W T=\Tr_{W_{\bar0}}T-\Tr_{W_{\bar1}}T$.
The supertrace of an odd endomorphism is zero, and supertrace extends linearly to all endomorphisms. The superdimension is
$\operatorname{sdim}W=\operatorname{str}_W\id_W$.

\subsection{Modules, duality and finiteness conditions}

An untwisted weak $V$-supermodule is a superspace $M$ with a
parity-preserving map
$a\mapsto Y_M(a,z)=\sum_{m\in\Z}a_m^Mz^{-m-1}$
satisfying truncation, the vacuum axiom and the module version of
\eqref{eq:super-jacobi}. All module morphisms are even. It is
\emph{admissible} if it has a grading $M=\bigoplus_{r\in\N}M(r)$,
compatible with parity, such that
$$
 a_m^M M(r)\subseteq M(r+\wt(a)-m-1)
$$
for homogeneous $a$; spaces with negative index are zero. The admissible
grading need not equal the $L(0)$-grading. Rationality means that every
admissible untwisted module is a direct sum of irreducible admissible
modules.

The \emph{parity reversal} $\Pi M$ has
$(\Pi M)_{\bar0}=M_{\bar1}$ and $(\Pi M)_{\bar1}=M_{\bar0}$, with
action $Y_{\Pi M}(a,z)\Pi m=\Pi(Y_M(a,z)m)$. This convention is
equivalent, by the parity operator, to the convention which includes
the factor $(-1)^{|a|}$ on the right. We call $V$ holomorphic if every
irreducible admissible untwisted supermodule is isomorphic to $V$ or $\Pi V$.
Thus we count irreducible modules up to parity reversal when discussing
holomorphicity.

An \emph{ordinary} module has semisimple $L(0)$, finite-dimensional
weight spaces and a lower-bounded weight grading. A grading-restricted
generalized module allows generalized $L(0)$-eigenspaces, still
finite-dimensional and lower bounded. Its contragredient $M'$ is the
direct sum of the duals of these generalized weight spaces. In the
categories used below, contragredience is an exact contravariant
equivalence, with $M''\cong M$. In particular, it preserves simplicity.
For an irreducible ordinary module over a strongly rational VOA we write
\begin{equation}\label{eq:ordinary-grading}
 M=\bigoplus_{r\in\N}M_{h_M+r},\qquad M_{h_M}\neq0.
\end{equation}
Here $h_M\in\Q$ is the lowest conformal weight. No positivity condition on $h_M$ is part of the definition of strong rationality.

An even bilinear form on $V$ is invariant if for homogeneous vectors
\begin{equation}\label{eq:invariant-form}
 \langle Y(a,z)u,v\rangle
 =(-1)^{|a||u|}
 \left\langle u,
 Y\bigl(e^{zL(1)}(-z^{-2})^{L(0)}a,z^{-1}\bigr)v
 \right\rangle.
\end{equation}
Such a form is supersymmetric. The VOSA is \emph{self-contragredient} if it admits a nondegenerate even invariant form, or equivalently if $V'\cong V$ by an even module isomorphism.
For a simple VOSA of CFT-type the invariant form is unique up to a
scalar when it exists, and self-contragredience gives $L(1)V_1=0$.
The normalization
$\langle\one,\one\rangle=-1$ gives
$u_1v=\langle u,v\rangle\one$ for $u,v\in V_1$, see
\cite[Section~4]{EM}.

For a VOA or VOSA $A$, and an $A$-module $M$, set
\begin{align*}
 C_2(A)&=\Span\{a_{-2}b:a,b\in A\},\\
 C_2^A(M)&=\Span\{a_{-2}m:a\in A,\ m\in M\}.
\end{align*}
The respective $C_2$-cofiniteness conditions mean that these subspaces
have finite codimension. The derivative identity implies
$a_{-r}m\in C_2^A(M)$ for $r\geq2$.

\begin{definition}\label{def:conventions}
An ordinary VOA is \emph{strongly rational} if it is simple,
self-contragredient, rational, $C_2$-cofinite and of CFT-type.
We use the same list of conditions for an integer-graded VOSA,
with rationality understood as above.
\end{definition}

In particular, neither unitarity nor positive lowest weights for
nonvacuum modules are included among these conditions.

\subsection{Vertex tensor categories}

We denote vertex tensor products by $\boxtimes$, and use $\otimes$
for ordinary vector-space or external tensor products.
For a simple, self-contragredient, $C_2$-cofinite VOA $A$ of CFT-type,
let $\Cat(A)$ be the category of grading-restricted generalized
$A$-modules. The logarithmic tensor theory and cofiniteness results give a braided abelian vertex tensor category with finite-length objects and finite-dimensional morphism spaces, closed under contragredients, see \cite{HuangFinite}. Its unit is $A$ and
its twist is
\begin{equation}\label{eq:twist}
 \theta_M=e^{2\pi iL(0)}.
\end{equation}
At this stage the category need not be semisimple or rigid.

Before establishing rigidity, we use the natural isomorphism
\begin{equation}\label{eq:GV-prelim}
 \Hom_{\Cat(A)}(P\boxtimes Q,A)
 \cong\Hom_{\Cat(A)}(P,Q').
\end{equation}
This is the Grothendieck--Verdier duality of
Allen--Lentner--Schweigert--Wood \cite[Theorem~2.12]{ALSW}, where
self-contragredience identifies the dualizing object $A'$ with $A$.
It gives the needed duality for intertwining operators without
asserting that every object already has a tensor dual.

If $A$ is strongly rational, Huang's theorem \cite{Huang} makes
$\Cat(A)=\Rep(A)$ a modular tensor category. Thus it is a semisimple
rigid ribbon category, has finitely many simple objects, and its
braiding is nondegenerate. Its simple objects are precisely the
irreducible ordinary $A$-modules. For a simple module $M$, the \emph{twist} is
$$
 \theta_M=e^{2\pi i h_M}\id_M.
$$
The categorical trace of $\id_M$ in this ribbon category is called the 
\emph{categorical dimension} of $M$, denoted by
$d_M$. Nondegeneracy means that the only simple object whose
monodromy with every object is the identity is the tensor unit.
The full subcategory of objects with this monodromy property is called the M\"uger center \cite{Muger}. Twist and braiding satisfy the balancing identity in modular tensor category.

An invertible simple object is called a \emph{simple current}.
Its tensor inverse is its contragredient, and the simple currents form
a finite abelian group under fusion. A category in which every simple
object is invertible is called \emph{pointed}. For a self-dual simple
current $J$, the self-braiding $c_{J,J}$ is a scalar on
$J\boxtimes J\cong A$. We use the \emph{spin--statistics identity}
\begin{equation}\label{eq:spin}
 c_{J,J}d_J=\theta_J,
\end{equation}
in the conventions of \cite[Corollary~2.8]{CKL}.
Note that categorical dimension and twist are different invariants.

\begin{lemma}\label{lem:invertible-sign}
If $g$ is invertible in a ribbon fusion category, then $d_g=\pm1$. 
\end{lemma}

\begin{proof}
Ribbon dimensions are invariant under duality and multiplicative, so $d_g=d_{g^*}$ and $d_gd_{g^*}=1$. This proves $d_g^2=1$.
\end{proof}

\section{\texorpdfstring{$C_2$}{C2}-cofiniteness of the even subalgebra}
\label{sec:c2}

Throughout this section let $V=U\oplus F$ be a simple, $C_2$-cofinite, integer-graded VOSA of CFT-type with $F\neq0$. Rationality and holomorphicity are not needed in this section.

\begin{lemma}\label{lem:inherited}
The even part $U$ is a simple VOA of CFT-type. The two summands $U$ and $F$ are nonisomorphic simple $U$-modules, and the coefficients of $Y(F,z)F$ span $U$. If $V$ is self-contragredient, then $U'\cong U$ and $F'\cong F$.
\end{lemma}

\begin{proof}
We apply McRae's \cite[Theorem~3.2]{McRorb} to the parity involution. The decomposition of $V$ gives simplicity of $U$, and irreducibility and inequivalence of the two parity summands $U$ and $F$ of $V$. The CFT-type condition passes to $U$ obviously.

Denote the coefficient span of $Y(F,z)F$ as $
I=\operatorname{span}
\{a_n b:a,b\in F,\ n\in\mathbb Z\}.
$
By the commutator formula, clearly $I$ is a $U$-submodule. If $I=0$, then $F$ is a nonzero proper ideal of $V$, contradicting the simplicity of $V$. Hence $I\neq0$, and the simplicity of $U$ implies $I=U$. 

When $V$ is self-contragredient, an even nondegenerate invariant form on $V$ restricts nondegenerately to its orthogonal parity summands. This gives
$U'\cong U$ and $F'\cong F$.
\end{proof}

The following cofiniteness statement is used in the proof of
\cite[Proposition~2.1]{CGNS}. It can be proved using Miyamoto's cofiniteness transfer theorem \cite{MiyTensor}.

\begin{proposition}\label{prop:c2}
The VOA $U$ is $C_2$-cofinite.
\end{proposition}

For completeness, we give the parity version of Miyamoto's argument. 
We first recall the $C_2$ case of \cite[Key Theorem]{MiyTensor}. Let $A$ be a VOA of CFT-type and let $M,N,W$ be $\N$-gradable $A$-modules. If $M$ and $N$ are $C_2$-cofinite and an intertwining operator of type $\binom{W}{M\,N}$ has coefficients spanning $W$, then $W$ is $C_2$-cofinite. This is the surjectivity for the intertwining operator, and no $C_2$-cofiniteness of $A$ is assumed. By Lemma~\ref{lem:inherited}, the coefficients of
$Y|_{F\otimes F}$ span $U$. Thus it suffices to prove that $F$ is
$C_2$-cofinite as a $U$-module.

Denote
$$
 C=C_2^U(F)=\Span\{u_{-2}a:u\in U,\ a\in F\},\qquad Q_F=F/C.
$$
Since $D(u_{-2}a)=(Du)_{-2}a+u_{-2}Da$, that is $DC\subseteq C$, the quotient $Q_F$ is a graded $\C[D]$-module, with $D$ of degree one.
Recall that a finitely generated $\C[D]$-module is finite-dimensional over $\C$ if and only if it is torsion. We will establish finite generation of $Q_F$ in Lemma~\ref{lem:D-finite} and torsion in
Lemma~\ref{lem:D-torsion}.

\begin{lemma}[{cf.\ \cite[Proposition~7]{MiyCyclic}}]\label{lem:D-finite}
There is a finite-dimensional homogeneous subspace $P\subset F$ such that
\begin{equation}\label{eq:D-finite}
 F=C+\C[D]P.
\end{equation}
\end{lemma}

\begin{proof}
For subspaces $A,B\subseteq V$, write $A_{-2}B=\Span\{a_{-2}b:a\in A,\ b\in B\}$. 
The equality $C_2(V)\cap F=U_{-2}F+F_{-2}U$ and $C_2$-cofiniteness of $V$
give a finite homogeneous space $P$ with $F=C+F_{-2}U+P$.
For homogeneous $a\in F$ and $u\in U$, skew symmetry gives
\begin{equation}\label{eq:minus-two-skew}
 a_{-2}u=-\sum_{j\geq0}\frac{(-1)^j}{j!}
                  D^j(u_{j-2}a).
\end{equation}
The $j=0$ term belongs to $C$. For $j>0$, the inner vector has weight $j$
less than $a_{-2}u$. Since $DC\subseteq C$, induction on weight proves
\eqref{eq:D-finite}. The grading is bounded below and the sum is finite by
truncation.
\end{proof}

We adapt the coefficient-function argument of \cite[Section~4]{MiyCyclic} to rule out a nonzero free summand of $Q_F$.

\begin{lemma}\label{lem:D-torsion}
The finitely generated $\C[D]$-module $Q_F$ is torsion. Consequently
$F/C_2^U(F)$ is finite-dimensional.
\end{lemma}

\begin{proof}
Suppose that $Q_F$ has a nonzero free part. The graded structure theorem over
$\C[D]$ gives homogeneous free generators. Choose one represented by
$\alpha\in F$, and quotient by the torsion part and the other free summands.
Write $\equiv$ for equality in this quotient. Its degree $\wt(\alpha)+x$ is
spanned by
\begin{equation}\label{eq:alpha-basis}
 \alpha_{-x-1}\one=\frac{D^x\alpha}{x!},\qquad x\in\N,
\end{equation}
and its lower degrees vanish.

For homogeneous $a,b$ of opposite parity, denote
$M=\wt(a)+\wt(b)-\wt(\alpha)$ and define $f_{a,b}:\N\to\C$ by
\begin{equation}\label{eq:coeff-f}
 a_{-x+M-1}b\equiv
 f_{a,b}(x)\alpha_{-x-1}\one.
\end{equation}
Let $\mathscr F_0$ and $\mathscr F_1$ be the spans of these functions for
$(a,b)\in U\times F$ and $(a,b)\in F\times U$, respectively. Translation gives
\begin{equation}\label{eq:Cx-stable}
 f_{Da,b}(x)=(x-M)f_{a,b}(x),
\end{equation}
so both spaces are $\C[x]$-modules. 
For $a\in U$ and $b\in F$, we have $a_{-r}b\in C$ whenever $r\ge2$. Hence the defining relation for $f_{a,b}$ gives
$f_{a,b}(x)=0$ for every $x\in\N$ with $x\ge M+1$.
Thus each generator of $\mathscr F_0$ has finite support,
and the same holds for every function in $\mathscr F_0$.

Define
$$
 (\mathcal Sf)(x)=\sum_{j=0}^{x}(-1)^j\binom{x}{j}f(j),
 \qquad (\mathcal Tf)(x)=(-1)^xf(x).
$$
Binomial inversion gives $\mathcal S^2=1$. Skew symmetry and
\eqref{eq:alpha-basis} give
\begin{equation}\label{eq:binomial-skew}
 f_{a,b}(x)=(-1)^M(\mathcal Sf_{b,a})(x).
\end{equation}
Here the lower degrees vanish and one input is even, so there is no super
sign. The binomial transform of a finitely supported function is polynomial,
hence
\begin{equation}\label{eq:F-polynomial}
 \mathscr F_1=\mathcal S\mathscr F_0,
 \qquad\mathscr F_1\subseteq\C[x]|_{\N}.
\end{equation}

We claim that
\begin{equation}\label{eq:T-stable}
 \mathcal T\mathscr F_1=\mathscr F_1.
\end{equation}
For homogeneous $u,v,\xi\in F$, $m\in\Z$ and $h\in\N$, put
$N=\wt(u)+\wt(v)+\wt(\xi)-\wt(\alpha)$.
The super Borcherds identity is
\begin{align*}
 \sum_{i\geq0}\binom{m}{i}(u_{r+i}v)_{m+n-i}\xi
 ={}&\sum_{j\geq0}(-1)^j\binom{r}{j}
       u_{m+r-j}v_{n+j}\xi\\
 &-(-1)^r(-1)^{|u||v|}
       \sum_{j\geq0}(-1)^j\binom{r}{j}
       v_{n+r-j}u_{m+j}\xi.
\end{align*}
Set $r=-x+N+h$ and $n=-m-2-h$. Since $u,v\in F$ are odd, $u_{r+i}v\in U$. Moreover,
$m+n-i=-2-h-i\leq-2$, so the translation identity gives
$(u_{r+i}v)_{m+n-i}\xi\in C_2^U(F)=C$.
Thus every summand on the left belongs to $C$, and the
left-hand side is zero in $Q_F$. 
This yields
\begin{align}\label{eq:borcherds-c2}
 0\equiv{}&\sum_{j\geq0}(-1)^j\binom{-x+N+h}{j}
 u_{-x+N+h+m-j}(v_{-m-2-h+j}\xi)\notag\\
 &+(-1)^{-x+N+h}\sum_{j\geq0}(-1)^j
 \binom{-x+N+h}{j}
 v_{-x+N-m-2-j}(u_{m+j}\xi).
\end{align}
%The plus sign before the second sum is the super sign.

Choose an integer $q>0$ such that $u_{m+j}\xi=0$ for $j\geq q$. For $0\leq h<q$, the
first sum has a cutoff independent of $x$ by truncation. Its coefficient
function lies in $\mathscr F_1$ by \eqref{eq:Cx-stable}. The coefficient matrix
of the $q$ alternating functions in the second sum is
\begin{equation}\label{eq:binomial-matrix}
 A_{h,j}(x)=(-1)^{N+h+j}\binom{-x+N+h}{j},
 \qquad 0\leq h,j<q.
\end{equation}
It is invertible over $\C[x]$, since
$$
 \det\left[\binom{t+h}{j}\right]_{0\leq h,j<q}=1.
$$
Here the Vandermonde determinant on $t,t+1,\ldots,t+q-1$ is
$\prod_{j=0}^{q-1}j!$, cancelling the column leading coefficients $1/j!$.
The signs give $\det A=(-1)^{qN}$, so $A^{-1}$ has polynomial entries.

Solving the coefficient equations obtained from
\eqref{eq:borcherds-c2} shows that the alternating coefficient
function $\mathcal T f_{v,u_m\xi}$ belongs to $\mathscr F_1$. 
The homogeneous coefficients $u_m\xi$ span $U$ by
Lemma~\ref{lem:inherited}. Taking linear combinations of the required weight
therefore gives $\mathcal T\mathscr F_1\subseteq\mathscr F_1$, and
$\mathcal T^2=1$ proves \eqref{eq:T-stable}.

But $f_{\alpha,\one}=1$, so \eqref{eq:F-polynomial} and
\eqref{eq:T-stable} would make $(-1)^x$ polynomial on $\N$. Its values at even
integers force that polynomial to be $1$, contradicting its values at odd
integers. Thus $Q_F$ is torsion, and finite generation makes it
finite-dimensional over $\C$.
\end{proof}

\begin{proof}[Proof of Proposition~\ref{prop:c2}]
By Lemma~\ref{lem:D-torsion}, $F$ is a $C_2$-cofinite $U$-module. Its integral,
lower-bounded grading is $\N$-gradable after a degree shift. The restriction
$$
 Y|_{F\otimes F}:F\otimes F\longrightarrow U((z))
$$
is a $U$-intertwining operator, since its Jacobi identity with an input from
$U$ has the ordinary even sign. Its coefficients span $U$ by
Lemma~\ref{lem:inherited}. 
Applying Miyamoto's cofiniteness transfer theorem recalled above
to this intertwining operator proves the proposition. 
\end{proof}

\begin{corollary}\label{cor:finite-U}
If $V$ is also self-contragredient, then $U$ is simple, self-contragredient,
$C_2$-cofinite, and of CFT-type. Thus $\Cat(U)$ is a braided abelian category
described in Section~\ref{sec:prelim}.
\end{corollary}

\section{Rationality of the even subalgebra}
\label{sec:sectors}
\label{sec:descent}
\label{sec:ribbon}

Assume the hypotheses of Theorem~\ref{thm:general} with $F\neq0$.
Corollary~\ref{cor:finite-U} gives the braided abelian category
$\Cat=\Cat(U)$ and the duality \eqref{eq:GV-prelim}. We prove
Proposition~\ref{prop:descent-intro} and determine the fusion group in this section. Rigidity of the whole category will be used only after rationality of $U$ has been established.

\begin{lemma}\label{lem:F-current}\label{lem:sectors}
The parity summands $U$ and $F$ satisfy
\begin{equation}\label{eq:F-data}
 F\boxtimes F\cong U,\qquad c_{F,F}=-\id,\qquad \theta_F=\id.
\end{equation}
There is an abelian direct-sum decomposition by monodromy with $F$,
\begin{equation}\label{eq:sector-decomp}
 \Cat=\Cat_+\oplus\Cat_-,
 \qquad
 \Cat_\epsilon\boxtimes\Cat_\delta
       \subseteq\Cat_{\epsilon\delta},
 \qquad U,F\in\Cat_+.
\end{equation}
In particular, there are no extensions between the two components.
\end{lemma}

\begin{proof}
The parity grading makes $V$ an abelian intertwining algebra with trivial
associativity cocycle and commutativity cocycle
$\Omega(a,b)=(-1)^{ab}$. The compact-orbifold tensor equivalence
\cite[Corollary~4.8]{McRorb} identifies the subcategory generated by its two parity summands with $\sVec$ as a braided tensor category. This gives $F\boxtimes F\cong U$ and self-braiding $c_{F,F}=-\id$. 
Since $V$ is integer-graded, $\theta_F=e^{2\pi iL(0)}=\mathrm{id}_F$.

We recall the monodromy decomposition associated with the
order-two simple current $F$, see \cite[Appendix~A.2]{CGNS}. 
Tensoring with $F$ is an equivalence. Thus there is a unique natural automorphism $\mu$ of the identity functor such that the monodromy
$$
 \M_{F,P}:=c_{P,F}c_{F,P}=\id_F\boxtimes\mu_P.
$$
Suppressing associativity isomorphisms, the hexagon identities and naturality of the braiding give
$$
 \M_{F,P\boxtimes Q}=\id_F\boxtimes\mu_P\boxtimes\mu_Q,
 \qquad
 \M_{F\boxtimes F,P}=\id_{F\boxtimes F}\boxtimes\mu_P^2.
$$
The second monodromy is the identity, since $F\boxtimes F\cong U$.
Consequently
$$
 \mu_{P\boxtimes Q}=\mu_P\boxtimes\mu_Q,\qquad \mu_P^2=\id_P.
$$
The natural idempotents $(\id_P\pm\mu_P)/2$ split in the abelian category,
giving \eqref{eq:sector-decomp}, including its tensor rule and the absence of
extensions between the components. Finally $\mu_U=\mu_F=1$, since the unit braids
trivially and $c_{F,F}^2=1$.
\end{proof}

\begin{lemma}\label{lem:neutral}
The category $\Cat_+$ is semisimple, with simple objects $U$ and $F$.
\end{lemma}

\begin{proof}
For $P\in\Cat_+$, local simple-current induction
\cite[Theorems~3.14--3.15]{CKL} gives an untwisted $V$-supermodule
\begin{equation}\label{eq:induction}
 \Ind(P)=V\boxtimes_U P,\qquad
 \Ind(P)|_U=P\oplus(F\boxtimes P).
\end{equation}
It is grading-restricted, lower bounded and of finite length. Shifting the finitely many conformal-weight cosets makes it admissible, so rationality and holomorphicity of $V$ decompose it into copies of $V$ and its parity reversal. The sum is finite because its restriction has finite length.
Hence \eqref{eq:induction} is a finite direct sum of copies of $U\oplus F$.
Its summand $P$ is therefore semisimple with constituents $U$ and $F$.
Both already belong to $\Cat_+$.
\end{proof}

\begin{proposition}\label{prop:klein}\label{prop:rank-four}
The VOA $U$ is strongly rational, and $\Rep(U)$ is pointed with fusion
group $\Z_2\times\Z_2$. Its odd simple current satisfies
\begin{equation}\label{eq:dF-negative}
 d_F=-1.
\end{equation}
\end{proposition}

\begin{proof}
Suppose first that $\Cat_-\neq0$ and choose a simple $X\in\Cat_-$.
Its contragredient $X'$ is also in $\Cat_-$, otherwise
Lemma~\ref{lem:neutral} and $U'\cong U$, $F'\cong F$ would place $X$ in
$\Cat_+$. Hence $X\boxtimes X'$ is a semisimple sum of $U$ and $F$.
The multiplicities are determined by \eqref{eq:GV-prelim}:
\begin{align*}
 \dim\Hom(X\boxtimes X',U)
   &=\dim\End(X)=1,\\
 \dim\Hom(X\boxtimes X',F)
   &=\dim\Hom(X, F\boxtimes X)=0.
\end{align*}
To justify the vanishing, write $\theta_X=t\id_X$ with $t\neq0$.
Balancing identity gives $\theta_{F\boxtimes X}=-t\id_{F\boxtimes X}$, so the simple objects $F\boxtimes X$ and $X$ are nonisomorphic.

It follows that $X\boxtimes X'\cong U$, and braiding gives the reverse product.
Thus $X$ is invertible, and tensoring with $X$ gives an equivalence $\Cat_+\simeq\Cat_-$, with inverse tensoring with $X'$. By Lemma~\ref{lem:neutral}, both components are semisimple. If $\Cat_-=0$, 
then $\Cat=\Cat_+$ is already semisimple. In either case $\Cat$ is semisimple,  so \cite[Proposition~4.16]{McRorb} implies rationality of $U$.
Together with Corollary~\ref{cor:finite-U}, this proves strong rationality.

Huang's theorem \cite{Huang} now makes $\Cat=\Rep(U)$ modular.
Nondegeneracy forces $\Cat_-\neq0$, since otherwise $F$ would be a
nontrivial simple object in the M\"uger center. The equivalence between the two components therefore
gives exactly four simple objects,
\begin{equation}\label{eq:four-simples}
 U,\quad F,\quad X,\quad F\boxtimes X,
\end{equation}
all invertible. The spin--statistics identity \eqref{eq:spin}, together
with \eqref{eq:F-data}, gives $d_F=-1$. Their fusion group cannot be
$\Z_4$, because a generator $g$ would satisfy $g^2=F$, but
Lemma~\ref{lem:invertible-sign} gives $d_F=d_g^2=1$.
It is therefore $\Z_2\times\Z_2$.
\end{proof}

This proves Proposition~\ref{prop:descent-intro}.

\section{Strict positivity of conformal weights of simple currents}\label{sec:positivity}

We prove Theorem~\ref{thm:positive-currents} for an arbitrary strongly
rational VOA $A$. For an irreducible $A$-module $X$, let $h_X$ denote its
lowest conformal weight, so that
$$
 X=\bigoplus_{n\geq0}X_{h_X+n},\qquad X_{h_X}\neq0.
$$
These weights are rational. Products and iterates of intertwining
operators converge in their usual domains and satisfy analytic
associativity and commutativity
\cite[Theorems~3.5 and~3.7, Remark~3.8]{HuangDE}.
Moreover, $\Rep(A)$ is a modular tensor category with
$\theta_X=e^{2\pi i h_X}\id_X$
\cite[Theorems~4.1, 4.5 and~4.6]{Huang}.
We use these results without any assumption on the signs of the
conformal weights of the other irreducible modules.

If $J$ is a simple current, write $JM=J\boxtimes M$. Its fusion inverse is the contragredient module $J'$, and therefore
\begin{equation}\label{eq:current-dual-weight}
 J^{-1}\cong J',\qquad h_{J^{-1}}=h_J.
\end{equation}
The proof compares the expansions of a four-point function at $0$, $1$ and $\infty$. We first need an elementary lemma in complex analysis.

\begin{lemma}\label{lem:polynomial-normalization}
Let $F$ be a nonzero multivalued holomorphic function on
$\C\setminus\{0,1\}$ whose analytic continuations span a
one-dimensional space. Suppose that compatible branches have convergent
expansions
\begin{align*}
 F(z)&=z^{\alpha_0}f_0(z) &&\text{near }0,\\
 F(z)&=(z-1)^{\alpha_1}f_1(z-1) &&\text{near }1,\\
 F(z)&=z^\gamma f_\infty(z^{-1}) &&\text{near }\infty,
\end{align*}
where $f_0,f_1,f_\infty$ are holomorphic at $0$. Then
$$
 P(z)=z^{-\alpha_0}(z-1)^{-\alpha_1}F(z)
$$
is a nonzero polynomial, and
\begin{equation}\label{eq:polynomial-degree}
 \gamma-\alpha_0-\alpha_1
 =\deg P+\operatorname{ord}_0(f_\infty)\in\Z_{\geq0}.
\end{equation}
\end{lemma}

\begin{proof}
The two scalar factors cancel the monodromies of $F$ around $0$ and $1$.
Since these loops generate the fundamental group of $\C\setminus\{0,1\}$, $P$ is single-valued.
The local expansions at $0$ and $1$ show that both singularities are
removable, so $P$ is entire. Put $D=\gamma-\alpha_0-\alpha_1$.
At infinity, up to a nonzero constant arising from the choice of branch,
$$
 P(z)=z^Dg(z^{-1}),\qquad
 g(w)=(1-w)^{-\alpha_1}f_\infty(w).
$$
If $s=\operatorname{ord}_0g=\operatorname{ord}_0f_\infty$, then
single-valuedness implies $D-s\in\Z$. If $D-s<0$, the entire function $P$ tends to zero at infinity. Thus $P$ is bounded on $\C$, Liouville's theorem makes it constant. 
Its limit at infinity then gives $P\equiv0$, contrary to the assumption that $F$ is nonzero.
Hence $D-s\geq0$. 
The same expansion gives polynomial growth of degree $D-s$, so Cauchy's estimates imply that $P$ is a polynomial of that degree.
\end{proof}

\begin{proposition}\label{prop:weight-inequality}
Let $J$ be a simple current of a strongly rational VOA $A$, and let $M$
be an irreducible $A$-module. Then
\begin{equation}\label{eq:current-weight-inequality}
 2h_J+2h_M-h_{JM}-h_{J^{-1}M}\in\Z_{\geq0}.
\end{equation}
\end{proposition}

\begin{proof}
Set $j=h_J$, $h=h_M$, $K=JM$ and $N=J^{-1}M$. We identify $J^{-1}$
with $J'$. By \cite[Proposition~4.1]{CKL}, there are lowest-weight
vectors $a\in J_j$, $b\in(J')_j$ and an intertwining operator
$\mathcal E$ of type $\binom{A}{J\ J'}$ such that, after rescaling,
\begin{equation}\label{eq:current-vacuum-coefficient}
 \mathcal E(a,t)b
 =t^{-2j}\bigg(\one+\sum_{n\geq1}v_nt^n\bigg),
 \qquad v_n\in A_n.
\end{equation}
Here self-contragredience identifies the target $A'$ in that proposition with $A$, and the CFT-type condition gives the leading term.

Choose $u\in M_h$ and $u'\in(M')_h$ with $\langle u',u\rangle=1$.
Consider the following matrix element of an iterate of intertwining operators, convergent for $0<|z-1|<1$,
\begin{equation}\label{eq:current-four-point}
 F(z)=\left\langle u',
       Y_M\bigl(\mathcal E(a,z-1)b,1\bigr)u\right\rangle.
\end{equation}
Its expansion at $1$ is
\begin{equation}\label{eq:current-expansion-one}
 F(z)=(z-1)^{-2j}f_1(z-1),\qquad f_1(0)=1.
\end{equation}
In particular, $F$ is nonzero.

The function extends analytically to the universal cover of
$\C\setminus\{0,1\}$, and analytic continuation around any loop
multiplies it by a scalar. To see this, the analytic associativity theorem identifies the spaces of products and iterates, whose common dimension is
$$
\dim\Hom_{\Rep(A)}\bigl(J\boxtimes(J'\boxtimes M),M\bigr)
=\dim\End_{\Rep(A)}(M)=1.
$$
In one product expansion, $N=J'\boxtimes M$ is the only intermediate module, and
each of the two intertwining-operator spaces is one-dimensional.
In the reversed order, the product expansion has $K=J\boxtimes M$ as its only intermediate module,
while the iterate has $A$ as its only intermediate module, the
corresponding intertwining-operator spaces are also one-dimensional.
Thus after choosing bases, the fusing and braiding isomorphisms
involved are represented by nonzero scalars
\cite[Section~2]{Huang}. These scalars depend only on the
intertwining operators, so evaluating on the chosen four vectors
gives a function with scalar monodromy. The argument requires no
assumption on the dimensions of the lowest-weight spaces.

In the domain $|z|>1$, associativity expresses $F$ up to a nonzero
constant as
$$
 \left\langle u',\mathcal Y_1(a,z)\mathcal Y_2(b,1)u\right\rangle,
$$
where $\mathcal Y_1$ and $\mathcal Y_2$ have types
$\binom{M}{J\ N}$ and $\binom{N}{J'\ M}$, respectively.
The component of $\mathcal Y_2(b,1)u$ in $N_{h_N+n}$ contributes a
power $z^{h-j-h_N-n}$. Convergence and the $L(0)$-grading consequently
give
\begin{equation}\label{eq:current-expansion-infinity}
 F(z)=z^{h-j-h_N}f_\infty(z^{-1}),
 \qquad f_\infty\text{ holomorphic at }0.
\end{equation}
After commutativity, the expansion in $0<|z|<1$ is similarly obtained
from
$$
 \left\langle u',\mathcal Y_3(b,1)\mathcal Y_4(a,z)u\right\rangle,
$$
with intertwining operators of types $\binom{M}{J'\ K}$ and
$\binom{K}{J\ M}$, respectively.
The component in $K_{h_K+n}$ now contributes
$z^{h_K+n-j-h}$. Hence
\begin{equation}\label{eq:current-expansion-zero}
 F(z)=z^{h_K-j-h}f_0(z),
 \qquad f_0\text{ holomorphic at }0.
\end{equation}
The intertwining operators are ordinary, so the expansions have integral
increments and no logarithmic terms. The constant terms of $f_0$ and
$f_\infty$ may vanish, which is allowed in
Lemma~\ref{lem:polynomial-normalization}.

Applying that lemma to
\eqref{eq:current-expansion-one},
\eqref{eq:current-expansion-infinity} and
\eqref{eq:current-expansion-zero}, with
$$
 \alpha_0=h_K-j-h,\qquad \alpha_1=-2j,
 \qquad \gamma=h-j-h_N,
$$
gives $\gamma-\alpha_0-\alpha_1=2j+2h-h_K-h_N\in\Z_{\geq0}$,
as required.
\end{proof}

\begin{proof}[Proof of Theorem~\ref{thm:positive-currents}]
Fusion with $J$ permutes the finitely many irreducible $A$-modules.
Fix an orbit of length $r$ and write
$$
 M_k=J^kM,\qquad h_k=h_{M_k},\qquad k\in\Z/r\Z.
$$
By Proposition~\ref{prop:weight-inequality},
$$
 d_k=2h_J+2h_k-h_{k+1}-h_{k-1}\in\Z_{\geq0}.
$$
Summing around the orbit yields
\begin{equation}\label{eq:current-orbit-sum}
 2r h_J=\sum_{k\in\Z/r\Z}d_k\geq0,
\end{equation}
so $h_J\geq0$.

Suppose that $h_J=0$. Every $d_k$ is then zero, and $h_{k+1}-h_k=h_k-h_{k-1}$. Periodicity makes this common difference zero. Since the orbit was arbitrary, $h_{JM}=h_M$ for every irreducible $M$. The balancing identity now gives
$$
 c_{M,J}c_{J,M}
 =\frac{e^{2\pi i h_{JM}}}{e^{2\pi i h_J}e^{2\pi i h_M}}
   \id_{J\boxtimes M}
 =\id_{J\boxtimes M}.
$$
By semisimplicity the same holds for every object of $\Rep(A)$, so $J$
is in the M\"uger center. However, nondegeneracy of $\Rep(A)$ implies that a 
simple object in M\"uger center is the tensor unit. Thus $J\cong A$. Conversely, $h_A=0$
by the CFT-type condition.
\end{proof}

From Theorem~\ref{thm:positive-currents}, we have the following useful  corollaries. 
\begin{corollary}\label{cor:pointed-positive}
If $A$ is strongly rational and $\Rep(A)$ is pointed, then
$h_M>0$ for every irreducible $A$-module $M\not\cong A$.
\end{corollary}

\begin{corollary}\label{cor:bound}
Let $J$ be a simple current of a strongly rational VOA $A$. If an
irreducible $A$-module has an orbit of length $r$ under fusion with $J$,
then $2r h_J\in\Z_{\geq0}$, and this integer is positive if $J\not\cong A$.
In particular, a non-vacuum simple current of fusion order $n\geq2$
satisfies
$$
 \begin{aligned}
  n h_J&\in\Z_{>0},&\qquad h_J&\geq\frac1n
      &&(n\text{ odd}),\\
  2n h_J&\in\Z_{>0},&\qquad h_J&\geq\frac1{2n}
      &&(n\text{ even}).
 \end{aligned}
$$
Both bounds are sharp for every indicated order.
\end{corollary}

\begin{proof}
The orbit assertion follows from \eqref{eq:current-orbit-sum} and
Theorem~\ref{thm:positive-currents}. For the vacuum orbit, put
$J^k=J^{\boxtimes k}$, $h_k=h_{J^k}$ and
$$
 d_k=2h_J+2h_k-h_{k+1}-h_{k-1},\qquad k\in\Z/n\Z.
$$
Proposition~\ref{prop:weight-inequality} gives $d_k\in\Z_{\geq0}$ and
$\sum_{k=0}^{n-1}d_k=2n h_J$. Since $J^{n-k}\cong(J^k)'$,
contragredience gives $h_{n-k}=h_k$. Thus $d_{n-k}=d_k$ and $d_0=0$.
For odd $n$, pairing the nonzero indices yields
$$
 n h_J=\sum_{k=1}^{(n-1)/2}d_k\in\Z_{>0},
$$
where strict positivity follows from
Theorem~\ref{thm:positive-currents}. For even $n$, the orbit sum already
gives $2n h_J\in\Z_{>0}$. The lower bounds follow.

To see sharpness, let $L=\Z\alpha$ with $(\alpha,\alpha)=d$, where
$d=2n$ for odd $n$ and $d=n$ for even $n$. In both cases $L$ is a
positive-definite even lattice and $\alpha/n\in L^*$. The strongly
rational lattice VOA $V_L$ has the simple current
$J=V_{L+\alpha/n}$ of exact fusion order $n$, by the fusion rule \cite[Corollary~12.10]{DL}. Its lowest conformal weight is
$$
 h_J=\min_{m\in\Z}\frac d2\left(m+\frac1n\right)^2
     =\frac{d}{2n^2},
$$
which attains the respective bounds.
\end{proof}

As a byproduct, Theorem~\ref{thm:positive-currents} also proves the positivity conjecture of M\"oller \cite[Conjecture~1.1]{Moeller}.
In the setting of the following corollary, M\"oller proved that the lowest
conformal weight of $V(g)$ is nonnegative and that its weight-zero
space is one-dimensional if the lowest conformal weight is zero
\cite[Theorem~1.4]{Moeller}. Here we are able to prove the strict positivity.

\begin{corollary}[M\"oller's positivity conjecture]\label{cor:moeller}
Let $V$ be a strongly rational holomorphic VOA and let
$g\neq\id$ be an automorphism of finite order. Then the unique
irreducible $g$-twisted $V$-module $V(g)$ has strictly positive lowest
conformal weight.
\end{corollary}

\begin{proof}
Let $n>1$ be the order of $g$ and put $A=V^{\langle g\rangle}$.
The cyclic-orbifold theorems \cite{MiyCyclic} and
\cite[Corollary~5.25]{CM} imply that $A$ is strongly rational.
By \cite[Theorem~5.2 and Proposition~5.6]{EMS}, its irreducible
modules are pairwise nonisomorphic simple currents $W^{(i,j)}$ 
with $i,j\in\Z/n\Z$, such that
$$
 W^{(0,0)}=A,\qquad
 V(g^i)\big|_A=\bigoplus_{j\in\Z/n\Z}W^{(i,j)}.
$$
These results do not assume positivity of the non-vacuum conformal
weights, see also \cite[Section~5]{Moeller}.
Since $n>1$, every $W^{(1,j)}$ is nonisomorphic to $A$.
Theorem~\ref{thm:positive-currents} gives $h_{W^{(1,j)}}>0$ for all $j$.
The conformal vector of $A$ is that of $V$, so
$$
 h_{V(g)}=\min_{j\in\Z/n\Z}h_{W^{(1,j)}}>0.
$$
\end{proof}

\section{Character asymptotics and the nonexistence theorem}
\label{sec:characters}

We complete the proof using the quantum-dimension formula of
Dong--Jiao--Xu \cite[Lemma~4.2 and equation~(4.1)]{DJX}, together with Huang's identification of categorical dimensions with ratios of entries of the VOA modular matrix.

Let $A$ be a strongly rational VOA, and denote its irreducible modules as $M_0=A,M_1,\ldots,M_s$. For $\tau$ in the upper half-plane, the character of $M_i$ is
\begin{equation}\label{eq:characters}
 \chi_i(\tau)=\Tr_{M_i}q^{L(0)-c/24},
 \qquad q=e^{2\pi i\tau}.
\end{equation}
These series converge there by Zhu's theorem
\cite[Theorems~5.3.2--5.3.3]{Zhu}. Let $S$ be the matrix of $\tau\mapsto-1/\tau$ on the indexed torus one-point trace functions. Specializing their insertion to the vacuum gives
\begin{equation}\label{eq:S-transform}
 \chi_i(-1/\tau)=\sum_{j=0}^s S_{ij}\chi_j(\tau).
\end{equation}
Huang's identification of the categorical and modular
$S$-matrices gives
\begin{equation}\label{eq:qdim-S}
 d_{M_i}=\frac{S_{i0}}{S_{00}},\qquad S_{00}\neq0,
\end{equation}
see \cite[Theorem~4.5]{Huang}. No hypothesis that nonvacuum
lowest weights are positive is used for this categorical identity.

If every nonvacuum irreducible $A$-module has strictly positive lowest
conformal weight, Dong--Jiao--Xu's result \cite[Lemma~4.2]{DJX} identifies the character-defined \emph{quantum dimension} with the same ratio:
\begin{equation}\label{eq:character-limit}
 \qdim_A(M_i):=\lim_{t\downarrow0}\frac{\chi_i(it)}{\chi_0(it)}
       =\frac{S_{i0}}{S_{00}}=d_{M_i}>0.
\end{equation}
The positive-weight hypothesis guarantees that the vacuum term dominates the modularly transformed characters. Section~\ref{sec:positivity} provides this hypothesis for the even subalgebra under consideration.

\begin{proof}[Proof of Theorem~\ref{thm:general}]
Suppose that $F=V_{\bar1}\neq0$. By
Proposition~\ref{prop:klein}, the even subalgebra $U$ is strongly
rational, its representation category is pointed, and $d_F=-1$.
Corollary~\ref{cor:pointed-positive} gives $h_M>0$ for every
irreducible $U$-module $M\not\cong U$. Therefore
\eqref{eq:character-limit} applies and gives
$$
 d_F=\lim_{t\downarrow0}\frac{\chi_F(it)}{\chi_U(it)}>0,
$$
contradicting $d_F=-1$. Hence $F=0$.
\end{proof}

\begin{remark}
The holomorphicity is essential to the proof of Theorem~\ref{thm:general}. Consider 
integer-graded affine VOSA $L_1(\mathfrak{osp}(1|2))$. It is strongly rational with two irreducible modules. Its even part is a strongly rational VOA isomorphic to $L_1(A_1)\otimes L(-3/5,0)$ with 8 irreducible modules. This even-part VOA is not pointed and has an  irreducible module with negative conformal weight. Thus our argument does not prevent the existence of such VOSA.
\end{remark}

\bigskip
\noindent
\textbf{Acknowledgements.}
KS thanks Thomas Creutzig for pointing the main question
of this paper to him in 2023, and Jethro van Ekeren and Ching Hung Lam for discussions. KS and HW thank the invitation of the conference “Vertex Algebras and Tensor Categories” at the Chern Institute of Mathematics in 2026 July, where the result was announced. 
KS is supported by China NSF grant No.~1250010384.

\bigskip
\noindent
\textbf{Declaration of AI use.}
OpenAI’s ChatGPT is used as an auxiliary tool for mathematical discussions, editorial revision and checks of references.

\end{document}